\def\acc{\rm acc }
\def\nacc{\rm nacc }

\def\Aun{\under {A}}
\def\Cun{\under {C}}
\def\proj{\rm pr}

\def\sqdi{\rlap{$\Box$}$\Diamond$}
\def\boxdi{\sqdi}
\documentstyle{jtbart}
\makeatletter
\def\@begintheorem#1#2{\it \trivlist \item[\hskip
\labelsep{\bf #2\ #1.}]}
\def\@opargbegintheorem#1#2#3{\it \trivlist
      \item[\hskip \labelsep{\bf #1\ #2\ (#3).}]}
\makeatother

\newtheorem{them}{Theorem}[section]
\newtheorem{quest}{Question}
\newtheorem{lemm}[them]{Lemma}
\newtheorem{cor}[them]{Corollary}
\newtheorem{prop}[them]{Proposition}
\newtheorem{thm}[them]{Theorem}

\newcommand{\jtbnumpar}[1]{\refstepcounter{them}
\trivlist
\item[\hskip \labelsep{\bf \thethm \ #1.}]}

\newcommand{\jtbdef}{\jtbnumpar{Definition}}
\def\jtbnot{\jtbnumpar{Notation}}
\let\endjtbdef=\par

\newtheorem{ssnote}{COMMENT}

\def\PROOF #1.{\par\noindent{\it Proof#1}.\ \ignorespaces}
\def\proof #1.{\par\noindent{\it Proof#1}.\ \ignorespaces}

\def\subm{\leq}

\def\Ascr{{\cal A}}

\let\?=\joinrel

\let\leftv=^
\let\rightv=^

\def\acc{\mathop{\rm acc}}
\def\nacc{\mathop{\rm nacc}}

\def\k/{\kern.2em}    

\def\dom{\mathop{\rm dom}}

\def\NF{\mathop{\rm NF}}
\def\otp{\mathop{\rm otp}}

\def\cpr{\mathop{\rm cpr}}

\def\cf{\mathop{\rm cf}}

\def\sup{\mathop{\rm sup}}
\def\min{\mathop{\rm min}}
\def\max{\mathop{\rm max}}

\def\Cscr{\cal C}
\def\down{\smash{\mathchar"0223}}

\let\union=\cup             %

\edef\bigcup{\mathop{\textstyle\mathchar\the\bigcup}}
\let\bigunion =\bigcup

\let\inter=\cap             %

\edef\bigcap{\mathop{\textstyle\mathchar\the\bigcap}}

\edef\bigwedge{\mathop{\textstyle\mathchar\the\bigwedge}}

\edef\bigvee{\mathop{\textstyle\mathchar\the\bigvee}}

\edef\sum{\mathop{\textstyle\mathchar\the\sum}}
\def\ind #1#2#3{#1 \mathbin{\mathop{\down}_{#2}} #3}

\def\math&{\ \& \ }

\let\imply\rightarrow        %
\def\force {\mathrel^\joinrel\rightarrow}
\def\force {\mathrel{\scriptstyle\mathrel^\joinrel\rightarrow}}
\def\forceq {\mathrel{\mathop{\force}\limits_{\textstyle\texsim}}}
\def\forceq{\mathrel^\joinrel
 \mathrel{\mathop{\rightarrow}\limits_{\smash{\textstyle\texsim}}}}
\def\forceq{\mathrel{\scriptstyle\mathrel^\joinrel
 \mathrel{\smash{\mathop{\rightarrow}\limits_{\smash{\raise
 2pt\hbox{$\scriptstyle\texsim$}}}}}}}
\let\all\forall             %

\let\exclaim=!                 %
\let\iso\approx             %
\let\texsim=\sim         %
\let\conj\sim             %
\def\conjp #1 {\conj_{#1}}     %
\let\sim\simeq             %
\let\neg=\lnot             %

\def\0bar{\bar 0}         %
\def\1bar{\bar 1}         %
\def\abar{\overline a}

\let\sat=\models                %

\def\bbar{\overline b}

\def\hbar{\overline h}

\def\vbar{\overline v}
\def\wbar{\overline w}
\def\xbar{\overline x}

\def\Mun{\underline M}
\def\Nun{\underline N}

\def\Bun{\underline B}
\def\Aun{\underline A}
\def\Cun{\underline C}

\def\Lbar{\bar L}

\def\Proof{Proof}

\def\Cscr{{\cal C}}

\def\Pun{\underline P}

\def\bigunion{\union}
\def\Mun{{\underline M}}
\def\Nun{{\underline N}}

\def\text#1{\ifmmode\leavevmode\hbox{#1}\else
   \typeout{Warning: \string\text \space used outside math mode!}
   \begingroup\hbox{#1}\endgroup\fi}

\title{Abstract Classes with few models have `homogeneous-universal'
models}
\author{J. Baldwin
\thanks{Partially supported by N.S.F. grant 90000139}
\\Department of Mathematics\\
University of Illinois, Chicago
\and
S. Shelah
Department of Mathematics\\
Hebrew University of Jerusalem
\thanks{Both authors thank the U.S. Israel
Binational Science foundation for its support of this project.
This is item 393 in Shelah's bibliography.}}
\begin{document}
\maketitle

This paper is concerned with a class {\bf K} of models and an abstract
notion of submodel $\subm$.  Experience in first order model theory has
shown the desirability of finding a `monster model' to serve as a
universal domain for {\bf K}.  In the original constructions of
J\'onsson and Fraiss\`e, {\bf K} was a universal class and
ordinary substructure played the role of
$\subm$.  Working with a cardinal $\lambda$ satisfying
$\lambda^{<\lambda} = \lambda$ guarantees appropriate downward
Lowenheim-Skolem
theorems; the existence and uniqueness of a homogeneous-universal
model appears to depend centrally on the amalgamation property.
We make this apparent dependence more precise in this paper.

The major innovation of this paper is the introduction of
weaker notion (chain homogeneous-universal) to replace the
natural notion of $({\bf K},\subm)$-homogeneous-universal model.
Modulo a weak extension of ZFC (provable if V=L),
we show (Corollary\ref{onehucor})
that a class {\bf K}
obeying certain minimal restrictions
satisfies a fundamental dichotomy:
For arbitrarily large $\lambda$,
either {\bf K} has the maximal number
of models in power $\lambda$ or {\bf K} has a unique chain
homogenous-universal model of power $\lambda$.
We show (\ref{uniquehom})
in a class with
amalgamation this dichotomy holds for the notion
of {\bf K}-homogeneous-universal model in the more
normal sense.

The methods here allow us to improve our earlier results
\cite{BaldwinShelahprimalii}
in two other
ways:  certain requirements on all chains of a given length are replaced
by requiring winning strategies in certain games; the notion of a
canonically prime model is avoided.
A full understanding of these extensions requires consideration of the
earlier papers but we summarise them quickly here.

Shelah emphasized in
\cite{Shelahuniversal}
that Tarski's union theorem has two components: closure  under
unions; each union is an amalgamation base (smoothness).  The first is
used to show the existence of a homogeneous universal model; the second
is needed for uniqueness.  In this paper we show that closure can be
replaced by the existence of a bound for each chain and even stronger
that
we need the boundedness only for a `dense' (in a sense made precise by a
game defined
below) set of chains.

In \cite{BaldwinShelahprimalii}
we established a dichotomy between the smoothness of a class and a
nonstructure theorem.  There was a weakness in our result;
although the definition of smooth (there is a unique  compatibility
class over each chain) does not involve the concept of a canonically
prime model we only established the theorem for classes equipped with a
notion of a canonically prime model.  We remedy that difficulty in this
paper at some cost.  First we
require
some additional set theoretic hypotheses (all provable if {\bf V=L}).
Second we must weaken the conclusion.  Instead of coding stationary
sets we can only guarantee that there are $2^{\lambda}$ models
of power $\lambda$.

The results here generalize an earlier result
proved in \cite{Shelahnonelemii}.
  That paper dealt with a class {\bf K} satisfying
the axioms discussed here but also closed under unions of {\bf K}-chains
and that was smooth.  Theorem 3.5 and Claim 3.4 of
\cite{Shelahnonelemii}
imply that if {\bf K}
is categorical in $\lambda$ and has few models of power $\lambda^+$ then
the unique model of power $\lambda$ is an amalgamation base.

We rely on many notations and definitions from
\cite{BaldwinShelahprimali} and \cite{BaldwinShelahprimalii}  but only
on rudimentary results from those papers.

Section~\ref{scene} contains the background notation.
In Section~\ref{games} of this paper we introduce several games; we are
able to express questions about the smoothness or boundedness of a class
{\bf K} in terms of winning strategies for these games.
Section~\ref{smbst} describes the set theoretic hypotheses necessary for
our construction.  We show in Section~\ref{manymodels} that a winning
strategy for Player NAM in Game 2 $(\lambda,\kappa)$ implies the
existence
of many models.  Section~\ref{mainthm} translates the
existence of a winning strategy for a player (B) trying to show
chains are bounded and the
failure of the player trying to force nonamalgamation
(NAM) to have a winning strategy into the existence and,
if there are few models, the uniqueness of
 `chain homogeneous-universal' models.
For a class with amalgamation this yields uniqueness and existence of
the $({\bf K},\subm)$ homogeneous-universal models.
In Section~\ref{concl} we summarise our results
and suggest some open problems.
Section~\ref{appendix} contains proofs of the combinatorial results
summarised in Section~\ref{smbst}.

We thank Chris
Laskowski and Bradd Hart for their valuable advice in preparing this
paper.
\section{Setting the Scene}
\label{scene}

Most of the notions used in this paper are defined in
\cite{BaldwinShelahprimali} or \cite{BaldwinShelahprimalii}.  They or
minor variants occur in earlier papers of Shelah, specifically
\cite{Shelahuniversal}.

{\bf $({\bf K},\subm)$ is an abstract class satisfying Axiom group A
of \cite{BaldwinShelahprimali}:}
\begin{description}
    \item[A0]  If $M \in K$ then $M \subm M$.
    \item[A1]  If $M \subm N$ then $M$ is a substructure of $N$
    \item[A2]  $\subm$ is transitive.
    \item[A3]  If $M_0 \subseteq M_1 \subseteq N$, $M_0 \subm N$ and $M_1
\subm N$ then $M_0 \subm M_1$.
\end{description}
We review here some of the less common
concepts.
All notions defined with cardinal parameters have the
obvious variants obtained by, e.g., replacing $\lambda$
by $<\lambda$.
{\bf K}$_{\lambda}$ is the class of members of {\bf K} with
cardinality $\lambda$.
A $(<\lambda,\kappa)$ chain is a {\bf K}-increasing  chain of
cofinality
$\kappa$
members of {\bf K}
($i<j$ implies $M_i \subm M_j$), each of cardinality
$<\lambda$.
A chain $\Mun$ is {\bf K}-bounded if there is an $M \in {\bf K}$
and a compatible family of maps $f_i$ mapping $M_i$ into $M$.
  {\bf K} is
$(<\lambda,\kappa)$-bounded if each
$(<\lambda,\kappa)$-chain is bounded.
  {\bf K} is
$(<\lambda,<\kappa)$-closed if the union of each such chain is in
{\bf K}.
Sections 2 and 4 of
\cite
{BaldwinShelahprimalii} contain a number of
examples that illustrate these concepts.

\jtbnumpar{
Assumptions}
\label{cardassump}
We fix for this paper a cardinal $\lambda$ with the following
properties.
\begin{enumerate}
\item $\lambda$ is a regular cardinal greater than the size of the
vocabulary of {\bf K}.
\item
There are no maximal
models in {\bf K}$_{<\lambda}$.
\item {\bf K} is $(<\lambda,\lambda)$-closed.
\item The $<\lambda$-L\"{o}wenheim Skolem property holds.
\item  For some regular $\kappa < \lambda$, some $S \subseteq \lambda$,
and cardinal $R\leq \lambda$,
$\text{\sqdi}_{\lambda,\kappa,R}(S)$ holds.
\end{enumerate}

We briefly justify this group of assumptions.
Property $\text{\sqdi}$ is explained in Section~\ref{smbst}.
If for some $\kappa$ and $\lambda$,
$\text{\sqdi}_{\lambda,\kappa,R}(S)$ holds then in particular
$\diamond(\lambda)$ holds and so $\lambda^{<\lambda} = \lambda$.
Thus, there will be no cardinality obstruction to obtaining
homogenous-universal models of power $\lambda$.
Assuming {\bf K} is $(<\lambda,\lambda)$ closed is virtually a
convention.  If we dropped that hypothesis we would reach the
conclusion that $2^{\lambda}$ structures of power $\lambda$ were
increasing unions
of members of {\bf K}$_{<\lambda}$.
We say {\bf K} has few models if there are $<2^{\lambda}$
isomorphism types of
models of
power $\lambda$.
We will need to assume a boundedness hypothesis
on {\bf K}; but it varies with particular theorems in this paper.
We fix the proper $\kappa, R$ and $\lambda$ in clause v) at the
appropriate time.

\section{Two Games}
\label{games}

This paper is crucially concerned with the question of exactly how one
bounds an increasing {\bf K}-chain of models.  We have learned that it
is
not essential to posit that all chains are bounded; only  a sufficient
number of them.  This sufficiency can be described in terms of winning
strategies for certain games.
Game 1 $(\lambda, \kappa)$
is played between player B who wants to establish that {\bf K}
is $(<\lambda,\leq \kappa)$-bounded and player NB who is trying
to show the opposite.

\jtbdef
\begin{enumerate}
\item
For an ordinal $\alpha$,
a play of Game 1 $(\lambda,\alpha)$
lasts at most
$\alpha$ moves.  During the play
Player B chooses models $\langle L_i: i <\delta\leq \alpha\rangle$;
Player NB chooses models $\langle P_i: i <\delta \leq \alpha\rangle$.
At move $\beta$,
\begin{enumerate}
\item Player B chooses a model $L_{\beta}$ in {\bf K}$_{<\lambda}$
that is a proper {\bf K}-extension of all the structures
$P_{\gamma}$ for $\gamma <\beta$.
\item Player NB chooses a model $P_{\beta}$ in {\bf K}$_{<\lambda}$
that is a proper {\bf K}-extension of $L_{\beta}$.
\end{enumerate}
Either player loses the game if at some stage
he does not have a legal move.
Player B wins if $\Pun$ is bounded; otherwise Player NB wins.
Player B has the advantage of playing first at limit ordinals;
the price is that he must guarantee the existence of a bound at each
limit stage.
\item We say Player B has a winning strategy for Game 1 $(\lambda,<R)$
if he has a uniform strategy to win all plays of Game 1 $(\lambda,\mu)$
for each ordinal $\mu < R$.
\end{enumerate}

\begin{prop}
If {\bf K} is $(<\lambda,\leq\kappa)$-bounded then Player B has
a winning strategy for Game 1 $(\lambda,\kappa)$.
\end{prop}


We continue to use the notations for properties of embeddings of chains
established in \cite{BaldwinShelahprimalii} and briefly reviewed in
Section~\ref{scene}.

\jtbdef
\begin{enumerate}
\item
 The {\bf K}-increasing
chain $\Nun$ of members of {\bf K}$_{<\lambda}$
extends $\Nun'$ if $\Nun'$ is an initial segment of $\Nun$
\item
Two extensions $\Nun'$ and $\Nun''$ of a chain $\Nun$
can be amalgamated over $\Nun$ if there is
an $N \in {\bf K}$ and embeddings of $\Nun'$ and $\Nun''$ into
$N$ which agree on $\Nun$.
\item
 The {\bf K}-increasing
chain $\Nun$ of members of {\bf K}$_{<\lambda}$
is a $\lambda$-{\em amalgamation base} if
$\Nun$ is bounded and
every pair of
{\bf K}-extensions of $\Nun$,
each in {\bf K}$_{<\lambda}$,
can be amalgamated over $\Nun$
in {\bf K}$_{<\lambda}$.
\end{enumerate}

To say $\Nun$ is an amalgamation base is to say {\bf K} is smooth at
the chain $\Nun$;
if {\bf K} satisfies the $<\lambda$-L\"owenheim Skolem property
the increasing
chain $\Nun$ of members of {\bf K}$_{<\lambda}$
is a $\lambda$-{\em amalgamation base} if and only if
all extensions of $\Nun$ are compatible over $\Nun$.

A canonically prime model over a chain  $\Mun$
can be viewed as a `strategy', choosing,
at each limit stage,
one
among a number of possible compatibility classes over an initial
segment of the
chain.
At stage $\delta$, the strategy depends not on the actual ordinal
$\delta$ but on the isomorphism type of $\Mun|\delta$.
Game 2 concerns the construction of chains that are amalgamation bases.
We rename Player B as Player NAM (he wants to prevent amalgamation) and
Player NB as Player AM (he wants to build an amalgamation base).

\jtbdef
\begin{enumerate}
\item
A play of Game 2 $(\lambda,\alpha)$
lasts $\leq\alpha$ moves.
The play of the game is exactly as in Game 1
(with NAM replacing B, AM replacing NB).
A player who has no legal move loses; however,
Player NAM wins instantly if
for some limit ordinal $\beta \leq \alpha$, $\Pun|\beta$ is bounded
but is
not a
$\lambda$-amalgamation base.
\item Game 2 $(\lambda,<\beta)$ is defined similarly, but Player NAM
wins only if for some $\alpha < \beta$, $\Pun |\alpha$ is bounded but
is not a $\lambda$-amalgamation base.
\end{enumerate}

Note that if player NAM wins instantly at some stage then the length of
the game is less than $\alpha$.

In both games
the decision about a completed
game is based only on the second
player's moves (NB or AM); apparently to win the first player must force
the second to make mistakes.  The following lemma shows this intuition
is misleading.

\begin{lemm}  Player B (NAM) wins Game 1 (Game 2) $(\lambda,\kappa)$ if
and only if  the sequence
$\Lbar$ constructed during the play
is bounded (and is not an amalgamation base).
\end{lemm}

\Proof.  A chain is bounded (an amalgamation base) if and only if each
cofinal subsequence is.

The preceding remark is quite straightforward; contrast it with the
difficulties involved in considering canonically prime models over
subsequences \cite{BaldwinShelahprimalii}.

Note that if player NAM plays a winning strategy for Game 2, he also
playing a winning strategy for Game 1.
\section{$\text{\sqdi}$}
\label{smbst}
We discuss in this section the principle $\text{\boxdi}$ from
\cite{Shelah221}. It is a
combination of
Jensen's combinatorial principles  $\Box$ and $\Diamond$
that will be used in our main construction.
Some justification is necessary for the use
of such strong set theory.  On the one hand, arguments with a strong set
theoretic hypothesis that conclude the existence of many nonisomorphic
models show
it is impossible to prove in ZFC that {\bf K} has few models.
Thus any `structure theory' that could be established for {\bf K}
in ZFC would have to allow the maximal number of models in a class
with `structure'.
On the other,
while for ease of statement we assert that these combinatorial
principles follows from V = L, in fact they only depend on the structure
of the subsets of $\lambda$.   Thus $2^{\lambda}$ can be as large
as desired while keeping \sqdi$_{\lambda,\kappa,R}$ for each
$\kappa < \lambda$
and some $R \leq \lambda$.
Moreover, the consistency of \sqdi\  can
be obtained by a forcing extension as well as by an inner model.
Some instances of \sqdi\ depend only on appropriate instances of GCH.
We begin by establishing some notation.

\jtbnot
\begin{enumerate}
\item  For any set of ordinals $C$, $\acc[C]$ denotes the set of
accumulation points of $C$, i.e., the $\delta \in C$ with $\delta = \sup
C \inter \delta$.  $\nacc[C]$ denotes $C-\acc[C]$.
\item
Fix regular cardinals
$\lambda >\kappa$; let $C^{\kappa}(S)$ denote the
set of $\delta \in S$ which have cofinality $\kappa$.  In the
following $\delta$ always denotes a limit ordinal.
\item  Suppose there is $S \subseteq \lambda$ and a collection
$\langle C_{\delta}:\delta \in S\rangle$.  Then for any
$S_1 \subseteq S$,
$\tilde S_1$ denotes $S_1 \union
\ \bigunion\{C_{\delta}:\delta \in S_1$\}.
\end{enumerate}
\jtbdef
[$\text{\sqdi}_{\lambda,\kappa,R}(S)$]
\label{boxdi}
We say that
$\Cscr = \langle
C_{\alpha}:\alpha \in S\rangle$
and
the sequence
$\Ascr = \langle A_{\alpha}:\alpha < \lambda \rangle$
witness that
$\lambda, \kappa, R$
satisfy
[$\text{\sqdi}_{\lambda,\kappa,R}(S)$]
if for some subset $S_1$ of $S$
the following conditions hold.
\begin{enumerate}
\item $\kappa$ is a regular cardinal $< \lambda$, $R$ is an ordinal
$\leq \lambda$.
\item  $S$ is stationary in $\lambda$ and contains all limit $\delta \in
\lambda$ with $\cf(\delta) < R$; $S$ contains only even ordinals.
\item $S_1 \subseteq C^{\kappa}(S)$.
\item  Each $C_{\delta} \subseteq S$.
 \item  If $\alpha \in S$
\begin{enumerate}
\item
$C_{\alpha}$ is a
closed subset of $\alpha$,
\item if $\beta \in C_{\alpha}$ then
$C_{\beta} = C_{\alpha} \inter \beta$,
\item $\otp(C_{\alpha}) \leq \max (R,\kappa)$, more precisely,
\begin{enumerate}
\item $\otp(C_{\alpha}) = \kappa$ if $\alpha \in S_1$ and
\item
$C_{\alpha} \inter \tilde S_1 =
\emptyset$ and
$\otp(C_{\alpha}) <R$ if $\alpha \in S- \tilde S_1$ is a limit ordinal,
\end{enumerate}
\item all nonaccumulation points of $C_{\alpha} $ are even successor
ordinals.

\end{enumerate}
\item  If $\delta \in S$ is a limit ordinal then
 $C_{\delta}$ is a club
in $\delta$.

\item
Each $A_{\alpha}\subseteq \alpha$ and
\begin{enumerate}
   \item if $\beta \in S_1$ and $\alpha \in \acc[C_{\beta}]$ then
$A_{\alpha} = A_{\beta}\inter \alpha$,
    \item for $A\subseteq \lambda$ and any closed unbounded subset $E$
of $\lambda$, $X_E = \{\delta \in S_1: \delta \in E \& \acc[C_{\delta}]
\subseteq E\ \&
\ A_{\delta} = A \inter \delta\}$ is stationary in
$\lambda$.
\end{enumerate}
\end{enumerate}

\jtbdef
$\text{\sqdi}_{\lambda,\kappa,R}$
holds if for some subset $S\subseteq \lambda$,
$\text{\sqdi}_{\lambda,\kappa,R}(S)$
holds.

We discuss the truth of this proposition for various choices of $\lambda$
and $\kappa$.  We begin with a case proved in ZFC; some cases of the
GCH are
needed to find cardinals satisfying the hypotheses.
Combining the methods of
\cite{Shelah221}, \cite{Shelahsquares},
and \cite{Shelah351} yields the following result; we include a full proof
in an appendix at
the suggestion of the referee.

\begin{lemm}
\label{justificationboxb}
Suppose  $\mu^{\kappa}= \mu$ and
$2^{\mu}=\mu^+ = \lambda$.
If $S^* \subseteq S^{\kappa}(\lambda)$ is stationary
then there is an $S \subseteq \lambda$ such that
$\text{\sqdi}_{\lambda,\kappa,\omega}$
holds and $S \inter C^{\kappa}(\lambda) \subseteq S^*$.
\end{lemm}

This combinatorial result is sufficient for our model theoretic
constructions if the problematical model theoretic situation
concerns chains of length $\omega$ ($\kappa = \omega$.  To deal with
chains of longer length (this is essential; see
\cite{BaldwinShelahprimalii}), the following stronger combinatorial
principal is needed.

\begin{lemm}[V=L]
\label{justificationboxc}
If $\lambda = \mu^+$ and $\kappa \leq \mu$, then
for some stationary $S$,
\boxdi$_{\lambda,\kappa,\lambda}(S)$.
\end{lemm}

\section{Players B and NAM construct many models}

\label{manymodels}

We show in this section that if for some
$\kappa < \lambda$ and an  ordinal
$R\leq \lambda$
Player B has a winning strategy for
Game 1 $(\lambda,\kappa)$ and
Game 1 $(\lambda,<R)$, {\bf K} is $(\lambda,\geq R)$-bounded,
and  Player NAM has a winning strategy for Game 2 $(<\lambda,\kappa)$
then
{\bf K} has the maximal number of models in power $\lambda$.
This is of course a technically weaker conclusion
than in
\cite{BaldwinShelahprimalii} where we showed that if {\bf K} is not
smooth then {\bf K} codes stationary sets.
But we have weaker hypotheses here and
this conclusion expresses a somewhat weaker intuition of nonstructure.
This weakening
of the result is reflected in a complication of the main argument.
In the earlier paper we constructed for each
stationary set $W$
a model $M^W$ coding $W$.
Here, we construct $2^{\lambda}$ models
simultaneously and destroy putative isomorphisms between them enroute.
The second author has in mind
a more elaborate version of our construction, which we don't expand on
here, that recovers the coding of stationary sets in this context.
The stronger conclusion in \cite{BaldwinShelahprimalii}
assumed the class {\bf K} was equipped with
notions of free amalgamation and
canonically prime models;  here we have no such assumption.

The many-models proof given here illustrates the role of canonically
prime model notion.  When $\cpr$ is in the formal metalanguage a
particular choice of limit model is specified so one can construct
a sequence of models and code a stationary set by asking the question,
`Is the limit model at
$\delta$ the canonically prime model?'  When we remove this notion from
the formal language we have to destroy isomorphisms between the possible
choices of a limit model.  Diamond allows us to do this.

We thank Bradd Hart for suggesting a simplification in the proof of the
main result.

\begin{thm}
\label{maintheorem:smooth3}
Fix regular cardinals $\kappa < \lambda$ and an ordinal
$R\leq \lambda$
satisfying $\text{\sqdi}_{\lambda,\kappa,R}$.
Suppose
\begin{enumerate}
\item
{\bf K} is $(<\lambda,\geq R)$-bounded;
Player B has a winning strategy for
Game 1 $(\lambda,<R)$.
\item
Player NAM has a winning strategy for
Game 2 $(\lambda,\kappa)$.
\item
Player B has a winning strategy for
Game 1 $(\lambda,\kappa)$.
\end{enumerate}
Then there are $2^{\lambda}$ members of {\bf K} with cardinality
$\lambda$.
Moreover, these models are mutually $\leq_{\bf K}$-non-embeddible.
\end{thm}

\jtbnumpar{Remark}  The connection between hypotheses i)
and iii) deserves some comment.
If $\kappa > R$ then iii) implies the second clause of i).
It is tempting to think that i) implies iii).  However, if $\kappa > R$
a game of length $> \kappa$ might still have cofinality $<R$.  Neither clause
of i) guarantees a winning stategy for that game.

\proof.
Fix
$S$,
$\langle C_i:i \in S\rangle$
and $\langle A_{\alpha}:\alpha < \lambda \rangle$ and $S_1$
to witness
$\text{\sqdi}_{\kappa,\lambda,R}$.

Using a pairing function and
condition vii) of Definition~\ref{boxdi} we can find $\langle
\nu_{\alpha},\eta_{\alpha},f_{\alpha}:\alpha < \lambda\rangle$ with
$\nu_{\alpha}, \eta_{\alpha}$ in $2^{\alpha}$ and $f_{\alpha}$ a
function
from $\alpha$ to $\alpha$ such that for any  $\nu,\eta \in\lambda$,
any function $f:\lambda \mapsto \lambda$, and any closed unbounded
subset $E$ of $\lambda$, for some $\alpha \in S_1 \inter E$,
$\acc[C_{\alpha}]\subseteq E$,
and for every $\beta \in \{\alpha\} \cup
\acc[C_{\alpha}]$ we have
$\nu_{\beta} = \nu \inter \beta$,
$\eta_{\beta} = \eta \inter \beta$, $f$ restricted to $\beta$ is
$f_{\beta}$.

By induction,
for each
$\alpha< \lambda$ and each $\nu \in 2^{ \alpha}$
we define a structure
$N_{\nu}$
whose universe is an ordinal $<\lambda$
so that if
$\alpha < \gamma$, $\nu \in 2^{\alpha}$, $\eta\in 2^{ \gamma}$
$N_{\nu}$ is a proper
{\bf K}-submodel of $N_{\eta}$.
Then we finish  the construction by defining for each
$\eta \in 2^{\lambda}$,
$N_{\eta}$ as $\union_{\alpha <\lambda} N_{\eta|\alpha}$.

When $\eta \in 2^{\alpha}$, we write
$\Nun_{\eta}$ for the sequence $\langle N_{\eta|\beta}:\beta <\alpha\}$.

We now define a winning strategy against the restriction of
a play of a game
to a club that is appropriate for the arguments here.

\jtbdef
Fix $\alpha \in S$.  Suppose $\Nun = \langle N_i:i < \alpha\rangle$
is a {\bf K}-increasing chain and $C$ is a club in $\alpha$ such that
$\nacc[C]$ is a set of even successor ordinals.
Consider an initial segment of
a play of Game 1 (Game 2) where $L_i =
N_{\gamma_i}$,  $P_i =
N_{\gamma_{i+1}-1}$, and $\langle \gamma_i:i< \alpha_0
\rangle$ is an enumeration
of $C$.  If $\alpha$ is a successor and $\alpha-1 \in C$,
say $C = \{\gamma_i:i \leq \alpha_0\}$ then $P_{\alpha_0}$ is not
defined by the preceding; let $P_{\alpha_0} = N_{\alpha -1}$.
(Otherwise, $P_i$ is well defined since by Definition~\ref{boxdi} v) d)
$\gamma_{i+1}$ is a successor ordinal.  Note $\gamma_{i+1}-1 \neq
\gamma_{i}$ because $C$ contains only even ordinals.)

Suppose that each
$L_i$ has been
chosen by Player B (NAM)'s winning strategy for this game and suppose
$N_{\alpha}$ is now chosen
according to Player B (NAM)'s winning strategy in this play of the game.
We say $N_{\alpha}$ has been chosen by playing Player B (NAM)'s winning
strategy for Game 1 (Game 2) on $\Nun|C$.

The notation $\Nun|C$ is slightly inaccurate since when $\alpha$ is a successor
the choice depends on $N_{\alpha-1}$.  Nevertheless we adopt the notation
because of its
suggestiveness in the limit ordinal case.

Note that if $\beta\in C_{\alpha}$
 then playing according to the winning strategy
on $\Nun|C_{\beta}$ is by condition of v) b of Definition~\ref{boxdi},
the same as playing a winning strategy on
an initial segment of $\Nun|C_{\alpha}$.

We use the following ad hoc notation.
We need to introduce this notation because we are using winning
strategies in both games 1 and 2 as the hypothesis for the construction.
If we simply assumed {\bf K} is $(\lambda,\lambda)$-bounded then we
wouldn't need this curlicue.

\jtbnot
For $\alpha < \beta \in S$,
$C_{\beta,\alpha}$ denotes $\{
\epsilon:
\epsilon \in C_{\beta}\,
 \&\, \epsilon \geq \alpha\}$.

\jtbnumpar{Construction}

We split the construction into several cases.  Let $$\tilde S_1 =S_1
\union \union
\{C_{\delta}:\delta \in S_1\}.$$

Most cases in the construction are defined by playing the winning
strategy of Player B or Player NAM on a closed unbounded subset of
$\beta$.
Condition v) of
$\text{\sqdi}_{\kappa,\lambda,R}(S)$ guarantees that the various cases
cohere.  Certain inductive properties of the construction are
incorporated in the description of the cases.

At stage $\beta$, we have  fixed  $\nu,\eta \in 2^{\beta}$
and a map $f_\beta$ from $\beta$ to $\beta$.  For each $\tau \in 2^{
\beta}$ we construct a model $N_{\tau}$.
In many cases $N_{\tau}$ is chosen by playing the winning strategy
for Player B or player NAM on $\Nun_{\tau}|C
_{\beta}$.  To see that these
strategies do not conflict note first that if every play on
$\Nun_{\tau}|C_{\beta}$ has been played by the winning strategy in
either game then it has been played according to a winning strategy for
Game 1 (as the winning strategy for Game 2 also wins game 1) and so
inductively $N_{\tau}$ can be chosen by Player B's winning strategy
in Game 1.  Moreover, since Player B has a uniform winning strategy
for all games of length less than $\lambda$ his play does not depend
on the particular game of length less than $\lambda$ that is being
considered.
The inductive hypothesis in the cases where Player NAM's
winning strategy in Game 2 is used are verified below.
The key step in the proof is subcase b) of Case IV.  The other stages
are preserving the induction hypothesis.

\begin{description}
\item[Case I.]
$\beta$ is a successor ordinal $\gamma +1$.  We will choose
$N_{\tau}$ as a proper
{\bf K}-extension of $N_{\tau|\beta}$.   This is possible since we
have assumed that there are no maximal models in {\bf K}$_{<\lambda}$.
In certain subcases however we must be more specific.
\begin{description}
\item[Subcase a.]  $\beta \in S- \tilde S_1$.
Choose $N_{\tau}$ by
playing Player B's
winning strategy for Game 1 on
$\Nun_{\tau}|C_{\beta}$.
\item[Subcase b.] $\beta \in \tilde S_1$.
If Player NAM has already won Game 2 played on
$\Nun_{\tau}|C_{\beta}$.
at some stage $\gamma < \beta$,
play the winning strategy for Player B on $C_{\beta,\gamma}$.  If not,
$\Nun_{\tau}|C_{\beta}$.
has been played by the winning strategy of Player NAM.
Choose $N_{\tau}$ by
playing the winning strategy for Player NAM for Game 2 on
$\Nun_{\tau}|C_{\beta}$.
 \end{description}
\item[Case
II.] $\beta$ is a limit ordinal and
$\beta \not \in S$.  By ii) of
Definition~\ref{boxdi},
$\beta$ has cofinality
at least $R$.
Choose $N_{\tau}$ to
bound $\Nun|\tau$
(by $(<\lambda,\geq R)$
boundedness) and with
$|N_{\tau}| < \lambda$ (by the $<\lambda$-L\"{o}wenheim Skolem
property).
\item[Case III.]
$\beta$ is a limit ordinal and
$\beta \in (S-\tilde S_1$). Since $\cf(C_{\beta}) < R$, we can
choose $N_{\tau}$
to bound $\Nun_{\tau}|\beta$,
by playing Player B's winning
strategy for Game 1
on
$\Nun_{\tau}|C_{\beta}$.
(Note $C_{\beta}$ is unbounded in $\beta$ by
Definition~\ref{boxdi}
vi).

\item[Case IV.]  $\beta \in \tilde S_1$; $\beta$ a limit ordinal.
The situation is interesting only if $\beta = \union \Nun_{\nu} = \union
\Nun_{\eta}$, $f_{\beta}|N_{\nu|\gamma}$ is a {\bf K}-embedding for each
$\gamma < \beta$ and $\tau = \nu$.
Unless all of these conditions hold, choose $N_{\tau}$ as in Case Ib)
by playing NAM's strategy on
$\Nun_{\tau}|C_{\beta}$ or
B's winning stategy on
$\Nun_{\tau}|C_{\beta,\gamma}$.
If they do hold there are two subcases.
 \begin{description}
     \item[Subcase a.]
$\Nun_{\tau}$ is
an amalgamation base.
If $\gamma \in C_{\beta}$ then by Subcase Ib, or Subcase IV b) at stage
$\gamma$, $N_{\tau|\gamma}$ was chosen by Player NAM's winning strategy
on $\Nun_{\tau}|C_{\gamma}$.  So, we are able to
apply Player NAM's winning strategy for Game 2
$(\lambda,\kappa)$ on $\Nun_{\tau}|C_{\beta}$ to choose $N_{\tau}$.
         \item[Subcase b.]
$\Nun_{\tau}$ is not
an amalgamation base.
$N_{\eta}$
has been defined.
Thus there are incompatible bounds $A_1, A_2\in {\bf K}_{<\lambda}
$ for $\Nun_{\tau}$.
If there is an extension $\hat f_{\beta}$
of $f_{\beta}$ such that
$\hat f_{\beta}(A_1)$ is compatible with $N_{\eta}$, $N_{\tau} = A_2$;
otherwise $N_{\tau} = A_1$.
     \end{description}

\end{description}

This completes the construction.
The various cases are clearly disjoint.
We finish the proof by proving the
following claim.

{\bf Claim.}  If $\sigma\neq \tau \in 2^{\lambda}$ then
there is no {\bf K}-embedding of
$N_{\sigma}$ into
$N_{\tau}$.

Suppose for contradiction that $f$ is such an embedding.
For $\alpha$ in
 a closed unbounded subset $C$ of $\lambda$,
$$\alpha =
\bigunion_{\beta <\alpha}N_{\sigma|\beta}
=\bigunion_{\beta <\alpha}N_{\tau|\beta}.$$
and and for each $\beta < \alpha$,
$f|N_{\sigma|\beta}$ is an {\bf K}-embedding into $N_{\tau|\beta'}$
for some $\beta' < \alpha$.
Thus, there is an $\alpha \in C \inter S_1$ with $f_{\alpha}=
f|{\alpha}$, $\nu_{\alpha} = \sigma|\alpha$, $\eta_{\alpha} = \tau
|\alpha$ and $\acc[C_{\alpha}] \subseteq C$.
Moreover, for $\beta\in \acc[C_{\alpha}]$, $f|\beta = f_{\beta}$.
Since $\alpha \in S_1$, $\cf(\alpha) = \kappa$.
Since Player NAM has
a winning strategy in Game 2 $(\lambda,\kappa)$,
either Player NAM wins Game 2 at stage $\alpha$
or Player NAM won at some earlier stage $\delta$.
In the first event,
 subcase IV b)
of the construction applied at stage $\alpha$,
guarantees there
is no embedding  of
$N_{\nu}=N_{\sigma|\alpha}$ into any extension of
$N_{\eta} = N_{\tau|\alpha}$.
In the second event,
$\delta\in \acc{C_{\alpha}}$  so by
condition vii) a of
Definition~\ref{boxdi} $A_{\delta} = A_{\alpha}\inter \delta$.  That is,
$f_{\delta}:\delta \mapsto \delta$ and since $f_{\alpha}$ is a
(sequence of)
{\bf K}-embedding(s) so is $f_{\delta}$.  Since $f_{\delta} = f|\delta$,
 subcase IV b)
applied at stage $\delta$
guarantees there is no extension of $f|\delta$ mapping
$N_{\sigma|\delta}$ into any extension of
$N_{\tau|\delta}$ and thus no extension of
$f|\alpha$ mapping
$N_{\sigma|\alpha}$ into any extension of
$N_{\tau|\alpha}$.
Thus there is no ${\bf K}$-embedding of
$N_{\sigma}$
into
$N_{\tau}$
and we finish.

\section{Classes with few models have `homogeneous-universal' models}
\label{mainthm}
We consider here several variants on the notion of homogeneous-universal
model and establish the existence and uniqueness of models satisfying
one of these notions for a class that has few models.
\jtbnumpar{Assumption}
\label{assumptions}
In this section we assume that the class {\bf K} has a notion of
strong submodel satisfying the  axioms of group A in
\cite{BaldwinShelahprimali} (listed in Section~\ref{scene}), the
properties of $\lambda$ enumerated in Section~\ref{scene} and
\begin{enumerate}
\item {\bf K} has fewer than $2^{\lambda}$ models of power $\lambda$.
\item {\bf K} is $(<\lambda,<\lambda)$-bounded.
\end{enumerate}
All the results of this section go through under these assumptions.
Since some of them require slightly less, many of the statements repeat
these overriding hypotheses or stipulate more technical conditions that
suffice.

The following obvious consequence of Theorem~\ref{maintheorem:smooth3}
is a key to this section.

\begin{lemm}
\label{notnam}
For a
regular cardinal $\kappa < \lambda$ such that some regular $R\leq
\lambda$, $\lambda$ satisfies $\text{\sqdi}_{\lambda,\kappa,R}$,
Assumption
~\ref{assumptions} implies that Player NAM does
not have a winning strategy for Game 2 $(\lambda,\kappa)$.
\end{lemm}

\jtbnumpar{Refinement}  Checking  the hypotheses for
Theorem~\ref{maintheorem:smooth3}, we see that instead of assuming
{\bf K} is {$(<\lambda,<\lambda)$-bounded} it suffices to assume {\bf K}
is {$(<\lambda,\geq R)$-bounded} and
Player B has a winning strategy for
Game 1 $(\lambda,<R)$.

Our key idea is to redo the Fraisse-Jonsson construction in a category
where models have been replaced by chains that are amalgamation bases.
While we can not derive the amalgamation property directly from a `few
models' hypothesis, we are able to derive the existence of a sufficient
number of amalgamation bases to carry out the construction.

\jtbdef  Let {\bf K}$^*_{<\lambda}$ be the class of
$(<\lambda,<\lambda)$-chains that are
amalgamation bases and such that every initial segment of the chain with
limit length is an amalgamation base.


We define a partial order on these chains.

\jtbdef Let $\Mun$ and $\Nun$ be $(<\lambda,<\lambda)$-chains.
\begin{enumerate}
\item
$\Mun \conj \Nun$ if $\bigcup \Mun = \bigcup \Nun$ and each
$M_i$ is {\bf K}-embedded in some $N_j$ by the identity and vice versa.
\item $\Mun \prec \Nun$
if $\Mun \conj \Nun$ or there is an $N_i$ in the sequence
$\Nun$ such that for each $j$ with $M_j \in \Mun$,
 $M_j  \subm  N_i$.  (In the latter case, we say
$\Nun$ properly extends $\Mun$.)
\end{enumerate}

\begin{lemm}
\label{ambaseexist}
\begin{enumerate}
\item
{\bf K}$^*_{<\lambda}$
is not empty.  Indeed, for each $N \in {\bf K}_{<\lambda}$,
 there is an $\Nun \in {\bf K}^*_{<\lambda} $ with $N_0=N $.
\item
For each $(<\lambda,<\lambda)$-chain
$\Mun$ there is an $\Nun \in {\bf
K}^*_{<\lambda} $
that properly extends $\Mun$.
\end{enumerate}
\end{lemm}

\Proof.
Consider any play of Game 2 $(\lambda,\kappa)$ where Player NAM chooses
$N$ as
$L_1$.  If Player NAM wins each
such play then he has a winning strategy for Game 2 $(\lambda,\kappa)$
played in the class of {\bf K}-extensions of $N$.  By
Lemma~\ref{notnam}, this contradicts
Assumption~\ref{assumptions}\,i).
If not, $\Nun =
\langle N,P_0,P_1, \ldots\rangle$
is the required member of ${\bf K}^*_{<\lambda} $.  For the second
claim, apply this argument to a bound for $\Nun$.

Now we can regard this collection of chains under this partial order,
$({\bf K}^*_{<\lambda},\prec)$,
analogously to our basic notion of an abstract class.  This class is
easily seen to have the amalgamation property.


\begin{lemm}  If
$\Mun^0 \prec \Mun^1$
and
$\Mun^0 \prec \Mun^2$ and each $\Mun^i$ is in {\bf K}$^*_{\lambda}$
then there is a model $N \in {\bf K} $
such  that both $\Mun^1$ and $\Mun^2$ can be embedded in $N$ over
$\Mun^0$.
\end{lemm}

\Proof.  Let $M^1$ be a bound of $\Mun^1$ and $M^2$ a bound for
$\Mun^2$.  Since $\Mun^0$ is a $\lambda$-amalgamation base
these two models and, a fortiori, the sequences can be amalgamated
over $\Mun^0$.

\medskip
Note that if $\Mun$ is {\bf K}-embedded in both $\Nun$ and $N$, then
$\Nun$ and $N$ can be amalgamated over $\Mun$.
Since each element of ${\bf K}^*_{<\lambda}$ is an amalgamation base
the
joint embedding property for members of
${\bf K}^*_{<\lambda}$
is an equivalence relation.  Formally

\jtbdef  Let $\Mun, \Nun \in
{\bf K}^*_{<\lambda}$.  Then
${\cal E}
(\Mun,
 \Nun)$ if there exists an $N$ such that both
$\Mun$ and $\Nun$ can be {\bf K}-embedded in $N$.

To avoid the notational inconvenience of dealing with compatibility
classes we posit:

\jtbnumpar{Assumption}
\label{chainassume}
${\bf K}^*_{<\lambda}$ has the joint embedding property.  I.e.,
For any $\Mun, \Nun \in
{\bf K}^*_{<\lambda}$,
there exists an $N$ such that both
$\Mun$ and $\Nun$ can be {\bf K}-embedded in $N$.

We consider four variants on the notion of homogeneous universal.  First
we list two variants on the normal notion; then we describe the
analogous versions for
${\bf K}^*_{<\lambda}$.
First we fix the meanings of universal.

\jtbdef\mbox{}
\begin{enumerate}
\item For any class {\bf K},
$M$ is ${\bf K}_{<\lambda}$-universal
if
all members of {\bf K} with cardinality less than $\lambda$ can be
{\bf K}-embedded in $M$.
\item
$M$ is ${\bf K}^*_{<\lambda}$-universal if
each member of
{\bf
K}$^*_{<\lambda}$ can be {\bf K}-embedded in $M$.
\item
$M$ is ${\bf K}_{<\lambda}${\em chain-universal}
 if
each $(<\lambda,<\lambda)$ chain
can be {\bf K}-embedded in $M$.
\end{enumerate}


Applying Theorem~\ref{ambaseexist}  and Assumption~\ref{assumptions} ii)
($(<\lambda,<\lambda)$-boundedness)
it is easy to see:

\begin{lemm}
\label{chuniv}
If
$M$ is ${\bf K}^*_{<\lambda}$-universal
then
$M$ is ${\bf K}_{<\lambda}$ chain-universal.
\end{lemm}


 \jtbdef\mbox{}
\label{homundef}
\begin{enumerate}
\item $M$ is
$({\bf K}_{\lambda},\subm)$-{\em homogenous-universal} if
$M$ is ${\bf K}_{<\lambda}$-universal
and for each
$N_0 \subm N_1\in {\bf K}_{<\lambda} $ any {\bf K}-embedding of $N_0$
into $M$ can be extended to a {\bf K}-embedding of $N_1$ into $M$.
\item $M$ is {\em strongly
$({\bf K_{\lambda}},\subm)$-homogenous-universal} if
$M$ is ${\bf K}_{<\lambda}$-universal
and each
isomorphism between {\bf K}-substructures of $M$ with power $<\lambda$
can be extended to an automorphism of $M$.
\item  $M$ is
{\em chain homogeneous-universal (for
${\bf K}^*_{<\lambda}$})
if $M$ is
${\bf K}^*_{<\lambda}$-universal and
for any pair $\Mun
\prec \Nun$ of members of ${\bf K}^*_{<\lambda} $, any {\bf K}-embedding
of $\Mun$ into $M$ can be extended to an embedding of $\Nun$ into $M$.
\item $M$ is {\em strongly chain homogeneous-universal (for
${\bf K}^*_{<\lambda}$)         }
 if
$M$ is ${\bf K}^*_{<\lambda}$-universal and
whenever the $(<\lambda,<\lambda)$-chains
$\Mun, \Nun \in {\bf K}^*_{\lambda}$
are isomorphic and {\bf K}-embedded in $M$, the isomorphism can be
extended to an automorphism of $M$.
\end{enumerate}

We quickly summarise the relations among these notions and then proceed
to supporting examples.  It is immediate from the definition that
a $({\bf K}_{\lambda},\subm)$-homogenous-universal model is
${\bf K}_{<\lambda^+}$-universal.   If the empty structure is
a {\bf K}-submodel of every structure the requirement of
${\bf K}_{<\lambda}$-universality in the definition of
a $({\bf K}_{\lambda},\subm)$-homogenous-universal model is
redundant.  Clearly, Definition~\ref{homundef} ii)
implies Definition~\ref{homundef}  i).
But in general Definition~\ref{homundef} ii)
is stronger than  Definition~\ref{homundef}  i).
If we do not require the counterexample to have
cardinality $\lambda$ this is true even in the first order case.  For
there are $\lambda$-saturated structures (thus satisfying i)) that are
rigid.  (E.g. rigid real closed fields \cite{Shelahsecondordiv}.)  These
notions
are
however equivalent if {\bf K} is smooth (and with
the $\lambda$-L\"{o}wenheim Skolem property).

Similarly iv) is stronger than iii).  We now show
notion i) is stronger than
notion  iii).

\begin{lemm}
If {\bf K} is $(<\lambda,<\lambda)$-bounded and
$M$ is a
$({\bf K_{\lambda}},\subm)$-homogenous-universal model then $M$ is
chain homogeneous-universal.
\end{lemm}

\Proof.
Let $\Mun$ be in
{\bf K}$^*_{<\lambda}$ and {\bf K}-embedded in $M$.  Suppose $\Mun \prec
\Nun$. By the
$<\lambda$-L\"{o}wenheim Skolem property there is an $M_1$ with
$|M_1|<\lambda$ and $\Mun$ is {\bf K}-embedded in $M_1 \subm M$.  By the
boundedness there exists a model $M_2$ that is a {\bf K}-extension of
$\Nun$.  Since $\Mun$ is an amalgamation base, $M_1$ and $M_2$ can be
amalgamated over $\Mun$ into some $M_3$ (with cardinality $<\lambda$ by
the $<\lambda$-L\"{o}wenheim Skolem property.)  Since $M$ is
$({\bf K}_{\lambda},\subm)$-{\em homogenous-universal}, there is an
embedding of $M_3$ into $M$ and the image of $\Nun$ verifies that
$M$ is chain homogeneous-universal.

\jtbnot  For some examples the $\lambda$ is not an important
parameter.  In these cases, we simply omit mentioning any parameter.
This means any reasonable choice of $\lambda$ works.

\jtbnumpar{Example}
We define several classes {\bf K} such that at
each $\lambda$,
{\bf K}
has
chain homogeneous-universal
models but has no
${\bf K}_{<\lambda}$-homogeneous-universal
models.
\begin{enumerate}
\item
Let {\bf K} be the class of all structures of the form $L_U =(L,<,R_U)$
that are described below
with
$\subm$ being the usual notion of substructure.
$(L,<)$ is a linear
order
with a last element.  For an arbitrary but fixed
subset $U$ of $L$, let $R_U(x,y)$  hold
if $x \in U\, \&\, x < y$. {\bf K} is the
class of structures $(L,<,R_U)$.
Each member of {\bf K} satisfies the following two sentences.

$$(\forall x) (\forall y) [R(x,y) \imply x <y].$$
$$(\forall x) [(\exists y) R(x,y) \imply (\forall z) [z > x \imply
R(x,z)]].$$

Now the key point is that if $M \in K$ and $x$ is not the last element
of $M$, $x \in U$ is a $\Delta^0_1$ property of $x$ so $U$ is preserved
by
extension except that the last element $e$ of $M$ may be in $U$ in some
extensions but not others.  To see $M$ is not an amalgamation base
choose an extension with $e\in U$ and another where it is not.

Now let $\Nun= \langle N_i:i<\omega\rangle$ be an increasing chain of
members of {\bf K} such that $N_{i+1}$ contains an element above all
those in $N_i$.
Then $\bigcup\Nun \not \in {\bf K}$ since it
doesn't have a last element.  But  $\Nun$ is in ${\bf K}^*$
since whether $a\in U$ has been determined for each $a\in \bigcup\Nun$.

Thus {\bf K} does not have any {\bf K}-homogenous-universal models (as
this implies the amalgamation property for models of smaller size.)  But
${\bf K}^*$ is nonempty and we obtain
${\bf K}^*$-homogeneous-universal models as in Theorem~\ref{huexist}.
This example has many models in each cardinality.

\item  Let $L$ have a single binary relation $R$.
Again let $\subm$ be the usual notion of substructure.
Let $T$ assert that
$R$ is asymmetric ($xRy \imply \neg yRx$) and irreflexive, that
each point is related to at most one other, and that every point but
one, which is related to no one, is related to exactly one other.  Then
$xRy \vee yRx \vee x=y$ defines an equivalence relation such that all
classes but one (with a single element) have two elements.  Now there is
a $\Delta^0_1$-definition of the predicate "$x$ is the `lower' of two
related elements" (provided $x$ is related to some element).  No model
in {\bf K} is an amalgamation base (the element in the 1-element class
can be `upper' or `lower' in an extension).  {\bf K} is
$(<\infty,<\infty)$-bounded but is not $(<\lambda,<\lambda)$-closed for
any $\lambda$.
For any $\lambda$, Player AM has a winning
strategy for Game 2 $(\lambda,<\lambda)$ so by the proof of
Theorem~\ref{huexist} there are
chain homogeneous-universal
models but as in the previous example they cannot be
${\bf K}_{<\lambda}$-homogeneous-universal.
{\bf K} is categorical in all infinite powers.
\item
If we vary this example by allowing {\bf K} to contain the structures
where all classes have two elements then {\bf K} is $(\infty,\infty)$
closed but the main point of the example continues to hold.  There are
now two models in each cardinality.
\end{enumerate}

\medskip
The last two examples show there is little hope of finding
${\bf K}_{\lambda}$-homogeneous-universal models in the abstract
setting solely from  hypotheses about the number of models.  We show
that such models exist and are unique under an amalgamation hypothesis
in Theorem~\ref{uniquehom}.

\jtbnumpar{Example}
Here is a class $({\bf K},\leq)$ and a model that is
${\bf K}_{<\lambda}$-homogeneous-universal
and
chain homogeneous-universal but {\bf K} is not smooth.
$M \in {\bf K} $ if $M$ is isomorphic to a substructure which is
closed under
subsequence of
$\langle\lambda^{\leq\omega+1},<,L_i\rangle$ where $<$ is the natural
order by
subsequence and the $L_i$ are predicates for levels.  $M\subseteq N$ is
a {\bf K}-substructure of $N$ if any infinite chain in $M$ that is
bounded in $N$ is also bounded in $M$.
Now if $M$ has
no finite branchs
and for each finite initial segment of a branch through
$M$ there are $\lambda$ branches with an upper bound and
$\lambda$ branches with no upper bound,
$M$ is
${\bf K}_{<\lambda}$-homogeneous-universal and so chain homogeneous
universal.  But any model (or chain of models) that contains an
infinite branch without an upper bound is not an amalgamation base.

For the remainder of this section we explore the properties of notions
iii) and iv).
We first prove the existence of
a ${\bf K}^*_{<\lambda}
$ chain homogeneous-universal model.
Then we will show using the assumption that
{\bf K} has
few models any such ${\bf K}^*_{<\lambda} $ chain homogenous-universal
model is strongly ${\bf
K}^*_{<\lambda}$ chain homogeneous and thus there is a unique ${\bf
K}^*_{<\lambda}$ chain homogeneous-universal model of power $\lambda$.



\begin{thm}
\label{huexist}
Suppose $\lambda$ satisfies the assumptions enumerated in
Section~\ref{scene}
and
that Player B has a winning strategy for
Game 1 $(\lambda,\lambda)$.
There is
a ${\bf K}^*_{<\lambda}
$ chain homogeneous-universal model with cardinality $\lambda$.
\end{thm}

\Proof.  We construct a sequence of models $\langle
M_i:i<\lambda\rangle$ of cardinality $< \lambda$
with the universe of $M_i$ contained in $\lambda$; the
required $M$ is $\bigcup M_i$.  Since $\lambda^{<\lambda}= \lambda$
(as
$\text{\sqdi}_{\lambda,\kappa,R}(S)$ holds),
we can let
$\langle
\Nun_i:i<\lambda\rangle$ enumerate a set of representatives (say with
universe contained in $\lambda$) of
all isomorphism types of
members of ${\bf K}^*_{<\lambda}$.

Having constructed $M_i$ we define $M_{i+1}$ with cardinality $<\lambda$
by a play of Game 1 $(<\lambda,\lambda)$.
We construct a sequence $\langle M_{i,j}: j < \alpha\rangle$ of
structures with universe contained in $\lambda$
where
$\alpha = |i| + |M_i|^{|M_i|} < \lambda$.
The $M_{i,j}$ will be the plays of Player B in the game.
The plays of his opponent, denoted $M'_{i,j}$ will guarantee the
homogeneity.
Let $G_j = \langle \Aun_j, \Bun_j, f_j\rangle$ for $j < \alpha$ be a
list of all triples such that $\Aun_j,\Bun_j \in {\bf K}^*_{<\lambda}$,
$\bigcup \Aun \subseteq \lambda$, $\bigcup \Bun \subseteq \lambda$,
$\Aun_j
\prec \Bun_j$, $\Bun_j \iso \Nun_{\beta}$ for some $\beta < i$,  and
$f_j$ maps $\Aun_j$ into $M_i$.  Now $M_i,0 = M_i$ and
$M'_{i,\delta} = M_{i,\delta}$ for $\delta$ a limit ordinal.
 $M'_{i,j+1}$ is
an amalgam of $M_{i,j}$ with $\Bun_j$ over $\Aun_j$ with the universe of
$M_{i,j+1}$ contained in $\lambda$.  $M_{i,j}$ is an extension of
$M'_{i,j}$ with universe a proper subset of $\lambda$ chosen by
Player B's winning strategy.  If $\delta$ is a
limit ordinal  $\leq\alpha$, $M_{i,\delta}$ is chosen as a
bound
for  $\langle M_{i,j}:j < \delta \rangle$ with cardinality $< \lambda$
by Player B's winning strategy.
$M_{i,\delta}$ of the proper cardinality can be chosen by the regularity
of $\lambda$ and the $<\lambda$-L\"{o}wenheim Skolem property.
Let
$M'_{i+1}$ be $M_{i,\alpha}$ and to
guarantee universality let $M_{i+1}$ be a common {\bf K}-extension of
$M'_{i+1}$ and $\Nun_{i+1}$.

The last result may be vacuous (e.g. nothing in the hypotheses
guarantees that {\bf K}$^*_{<\lambda}$ is not empty.)  Adding
Assumption~\ref{assumptions} and applying Lemma~\ref{chuniv} we deduce:

\begin{cor}
\label{huexistcor}
Suppose $\lambda$ and {\bf K} satisfy the assumptions enumerated in
Paragraphs \ref{cardassump}, \ref{assumptions},
\ref{chainassume} and
$\lambda^{<\lambda} =
\lambda$.  There is
a ${\bf K}^*_{<\lambda}
$ chain homogeneous-universal model with cardinality $\lambda$
that is {\bf K}$_{\lambda}$ universal.
\end{cor}

We will now show that if {\bf K} has few models of power $\lambda$ then
any
chain homogeneous-universal model
is actually
strongly chain homogeneous-universal
and thus the uniqueness of the
chain homogeneous-universal model.
The keys to showing this
lemma are the following invariants.

Expand
the language $L$ of our class $K$ to $L^*$ by adding for each $\alpha <
\lambda$ an $\alpha$-ary predicate $P_{\alpha}$.  Expand each
$L$-structure $M$ to an $L^*$-structure $M^*$ by letting $P_{\alpha}$
hold of an $\alpha$-sequence $\abar$ just if $\abar$ enumerates a
{\bf K}-submodel of $M$.  Without loss of generality the universe of $M$
is $\lambda$.
As in \cite{Shelah129}, for each $\beta <
\lambda$, $\phi_{\beta^*}$ describes the $\infty\lambda$-type of
$\beta$  (i.e. of $\langle \gamma:\gamma < \beta\rangle$).
We denote by $\xbar_j$ the sequence of variables
$\langle x_i: i<j\rangle$.
\jtbdef
For any representation $\Mun$ of a model $M$,
$$S(M,\Mun) =
\{i:i < \lambda \text{ is a limit ordinal and }  M \sat \neg (\all
\xbar_i)
\bigwedge_{j<i} \phi_{j^*}(\xbar_j) \rightarrow \phi_{i^*}(\xbar_i)
)\}.$$
$$\hat S(M,\Mun) =
\{i:i < \lambda \text{ is a limit ordinal and }  \langle M_j:j <
i\rangle  \not
\in {\bf K}^*_{<\lambda}\}.$$

It is easy to see that $S(M,\Mun)$ is an invariant of $M$ modulo the
cub
filter; thus we abbreviate $S(M,\Mun)$ to $S_M$.  The same holds for
$\hat S$.  Thus, $S_M$ is exactly the invariant defined in
\cite{Shelah129}.
(It appears to be the complement but this is a typo in the
earlier paper.)  For $X \subset \lambda$, we write $X^c$ for
$\lambda -X$.

\jtbnumpar{Claim}  If $M$ is chain homogeneous-universal $S_M$ and
$\hat S_M$ are equivalent modulo the cub filter on $\lambda$.

\Proof.
We say that an $L_{\infty,\lambda}$
formula $\phi(\xbar)$ is complete if it implies the quantifier free
diagram of $\xbar$; denote this formula by $\phi^*(\xbar)$.
Now we will show that
for each complete formula $\phi(\xbar) \in L_{\infty,\lambda}$ such that
$\phi(\xbar)$ implies that $\xbar$ enumerates an {\bf K}-increasing
chain,
and for
any sequence $\abar \in M$ that enumerates a member of {\bf
K}$^*_{<\lambda}$,
$$M \sat \phi(\abar) \text{ if and only if } M \sat \phi^*(\abar).$$

We prove this assertion by induction on quantifier rank; for quantifier
free formulas it is tautologous.
Suppose $\phi(\xbar)$ is $(\exists \vbar) \theta(\vbar,\xbar)$ and by
induction that we have the proposition for $\theta$.  Fix any $\bbar
\in M$ that enumerates a member of {\bf K}$^*_{<\lambda}$ and $M \sat
\phi^*(\bbar)$.  We must show $M\sat \phi(\bbar)$.  Suppose
$\phi(\abar)$ is witnessed by $\theta(\abar_1,\abar)$.
By Lemma~\ref{ambaseexist} there is an $\abar_2 \subseteq M$
such that
$\abar_1 \abar_2 \abar$
enumerates a member of ${\bf K^*_{<\lambda}}$.
Let
$\delta (\vbar,\wbar,\xbar)$
be the conjunction of the
quantifier-free diagram of
$\abar_1 \abar_2 \abar$
and $\phi'(\vbar,\wbar,\xbar)$  denote $\theta(\vbar,\xbar) \wedge
\delta (\vbar,\wbar,\xbar)$.
Since $M$ is chain homogeneous-universal, there exist
$\bbar_1,\bbar_2 \subseteq M$
such that
$\bbar_1 \bbar_2 \bbar$
enumerates a member of ${\bf K^*_{<\lambda}}$ isomorphic to
$\abar_1 \abar_2 \abar$.  By induction, we have $\sat
\phi'(\bbar_1,\bbar_2,\bbar)$ and thus $\sat \phi(\bbar)$ as required.

Now suppose $i\in \hat S_M^c$.  Then $M|i$ is an amalgamation
base and we have just observed that the
$\infty \lambda$ type of $i$ is equivalent to a quantifier free formula.
So $i$ is certainly not in $S_M$.
Thus $\hat S_M^c \subseteq S^c_M$.

The set $C = \{\delta < \lambda: \delta = \sup (\delta \inter
\hat S^c_M) \& \delta \text{ is the  universe of } M|\delta\}$ is a cub.
To complete the proof, it suffices to show $\hat S_M \inter C \subseteq
S_M$.  (This yields $\hat S_M = S_M$ modulo the cub filter.)
But if $i \in \hat S_M$,
$\Mun|i \not\in {\bf K}^*_{<\lambda}$.  So there exist
$(<\lambda,<\lambda)$ chains $\Mun'$ and $\Mun''$ whose restriction to
$i$ are isomorphic to $\Mun$ but which are incompatible over $\Mun|i$.
Since $M$ is $<\lambda$ chain-universal there are copies of both $\Mun'$
and $\Mun''$ {\bf K}-embedded in $M$.
Note
that if $i\in C$,
the $\infty \lambda$ type of $i$ is equivalent to a quantifier free
formula.
Thus $\Mun'|i$ and $\Mun''|i$ have
the same $\infty \lambda$-type and so
witness that $i  \in
S_M$.

\jtbnumpar{Assumption}  We assume for the next theorem the combination
of $\Box$ and $\Diamond$ described on page 7 of \cite{Shelah129}.

\begin{thm}  If there is an $M\in {\bf K} $ such that $\hat S_M$ is
stationary
then there are $2^{\lambda}$ models of {\bf K} with power $\lambda$.
\end{thm}

\Proof.  Since $\hat S_M = S_M$
we can just quote the main result of \cite{Shelah129}.  It
is shown in that paper, that for $\lambda$ satisfying the set theoretic
principle described there,
each model $M$
of power $\lambda$ is categorical in the language
$L_{\infty,\lambda}$ if $S_M$ is not stationary
or has $2^{\lambda}$ models of power $\lambda$ that
are $L_{\infty,\lambda}$-equivalent to it if $S_M$ is stationary.
The $2^{\lambda}$ models defined in
the proof are all unions of members of ${\bf K}_{<\lambda}$ so since
{\bf K} is $(<\lambda,\lambda)$-closed we have the result.

\begin{thm}  Suppose $M$ has cardinality $\lambda$, $\hat S_M$
is not stationary, $M$ is chain homogeneous-universal then
$M$ is strongly chain homogeneous-universal.  More generally, if
$M$ and
$N$ are any two chain homogenous-universal models and neither $\hat S_M$
nor $\hat S_N$ is stationary then $M$ and $N$
are isomorphic.
\end{thm}

\Proof.  This is the other half of the dichotomy proved in
\cite{Shelah129}.  The argument is presented in detail there.   To
summarise, suppose $M$ and $N$ are each ${\bf K}^*_{<\lambda}
$ chain homogeneous-universal.  Construct by induction a back and forth
between $M$ and $N$.  The successor stages are easy by the definition of
${\bf K}^*_{<\lambda} $-homogeneous.  The tricky point is the limit
stage.  But since $\hat S_M$ and $\hat S_N$ are not stationary we can
(by
restricting to a cub) assume they are empty.  Thus at limit stages the
sequences already constructed in $M$ and $N$ are in ${\bf
K}^*_{<\lambda} $.  So we can apply the definition of ${\bf
K}^*_{<\lambda} $-homogeneity and continue the construction.

\begin{cor}
\label{onehucor}
Suppose $\lambda$ and {\bf K} satisfy the assumptions enumerated in
Paragraphs \ref{cardassump}, \ref{assumptions},
\ref{chainassume}.
There is
a unique ${\bf K}^*_{<\lambda}
$
chain homogeneous-universal
model with cardinality $\lambda$
and it is
${\bf
K}^*_{<\lambda} $
strongly chain homogeneous-universal.
\end{cor}

\begin{thm}
\label{uniquehom}
Suppose for each $\kappa < \lambda$ and some $R\leq \lambda$, {\bf K}
and $\lambda$
satisfy the assumptions
enumerated in Paragraphs \ref{cardassump} and \ref{assumptions}.
Suppose further that
{\bf K}
satisfies the amalgamation and joint
embedding propeties and
$\lambda^{<\lambda} =
\lambda$.  
There is
a unique $({\bf K}_{\lambda},\subm)
$
homogeneous-universal
model with cardinality $\lambda$
\end{thm}

\Proof.  With amalgamation it is straightforward to construct
an $({\bf K}_{<\lambda},\subm)$ homogeneous-universal model.
It is chain homogeneous-universal and thus unique.

While we assumed \sqdi$_{\lambda,\kappa,R}$ for each $\kappa < \lambda$,
this is excessive.  We only need this assumption for those $\kappa$
where there is a possible failure of smoothness.
\medskip
The following result does not use the hypothesis that {\bf K} has few
models.
\begin{thm}  The following are equivalent.
\begin{enumerate}
\item $S_M$ is not stationary for some
${\bf
K}^*_{<\lambda} $ chain
homogeneous-universal model $M$ of power $\lambda$;
\item
Player AM has a winning strategy in Game 2 $(\lambda,<\lambda)$.
\end{enumerate}
\end{thm}

\Proof. First we show i) implies ii).
Let $\Mun$ be a
representation of $\Mun$ such that $S_M$ is not stationary; fix a closed
unbounded set $C$ that is disjoint from $S_M$.
To win player AM chooses $P_i$ as
a {\bf K}-extension
of $L_i$ and an isomorphism $f_i$ (extending the $f_j$ for $j<i$)
of $P_i$ with
an $M_{j_i}$
whose universe is an ordinal in $C$.
Thus, each limit stage in the chain
constructed  by the play of the game is in $C$ and so
is an amalgamation
base as required.

To see ii) implies i) modify the proof of Theorem~\ref{huexist}.
Construct a sequence $M'_{i,j}$ as follows.  Play Game 2
$(\lambda,<\lambda)$.  Let player NAM's moves be $L_j = M_{i,j}$ where
$M_{i,j}$ is chosen just as in the proof of Theorem~\ref{huexist}.  Let
$P_j = M'_{i,j}$ be chosen according to Player AM's winning strategy.
This guarantees that letting  $M' = \bigcup_{i,j} M'_{i,j}$,
$S(M',\Mun')$ is not stationary (indeed empty).

\jtbnumpar{Extension}  If $\lambda = \kappa^+$ and $\Box_{\lambda}$
holds ($\Box_{\kappa}$ in Jensen's notation) then Condition ii) in the
last lemma can be replaced by Game 2 $(\lambda,\kappa)$.
\section{Conclusions and Problems}
\label{concl}

Much of this paper can be viewed as answering the question:  What is the
role of closure under union in the construction of homogeneous-universal
models?  We have the following symbolic equation:
$$\text{closure under unions} = \text{boundedness} +
\text{smoothness}.$$
We have shown that the boundedness hypothesis can be weakened to the
existence of a winning strategy for Player B in Game 1; this suffices
to show the existence of
{\bf K}$^*_{\lambda}$
chain homogeneous-universal
models.  Similarly, we have weakened smoothness to the nonexistence of a
winning strategy for Player NAM in game 2; this suffices to prove the
chain homogeneous-universal
is ${\bf K}_{<\lambda}$ chain-universal.
We are then able to apply the argument from
\cite{Shelah129} to show that if {\bf K} has few models the chain
homogeneous-universal model is unique.  The most obvious question is

\begin{quest}  Can the results of this paper be obtained in ZFC?
\end{quest}

Rami Grossberg has made progress on this question and several related
ones.  In particular he has shown that there is a model of ZFC + CH +
$2^{\aleph_1} = 2^{\aleph_2}$ in which there is a class {\bf K} that
satisfies the model theoretic conditions of
Theorem~\ref{maintheorem:smooth3} (with $\lambda = \aleph_2$)
but has only one model of power
$\aleph_2$.  In the case where $\lambda$ is a successor cardinal
Grossberg has weakened our set theoretic assumptions.  These
results are still being written up.

But there is another use of closure under unions.  Separating one of
the components of the notions of limit model defined in
\cite{Makowskyabstractembed,Shelahnonelemii} and generalizing
\cite{AlbertGrossbergrich}, call a structure $M$ {\bf K}-rich if it
is {\bf K}-isomorphic to a proper substructure.  Now it is
easy to see that if a structure of power $\lambda$ is rich and ${\bf
K}_{\leq \lambda}$
universal and if {\bf K} is closed under union then there is a member of
{\bf K} with power $\lambda^+$.  There are examples \cite[page
431]{Shelahnonelemii} of classes {\bf K} that are $\lambda$-categorical
such that the model of power $\lambda$ is universal, maximal, and
homogenous (since rigid).

\begin{quest}  Suppose {\bf K} has few models in power $\lambda$.
What are minimal model theoretic conditions on a class
{\bf K} so that the homogeneous-universal model is not maximal?
Can the homogeneous-universal model be a Jonsson model?
\end{quest}

\begin{quest}
What sort of transfer theorems can be proved for the existence of
{\bf K}$^*_{\lambda}$-homogeneous-universal models (in various $\lambda$).
\end{quest}

Much of our efforts have been dedicated to showing the uniqueness of the
homogeneous-universal model.  This raises a metatheoretical question.

\begin{quest}  How much of stability theory can be carried through in
abstract class {\bf K} that has homogeneous-universal models (in many
cardinalities) but where they may not be unique.
\end{quest}

\section{Appendix:  Set Theoretic Lemmas}
\label{appendix}
In this section we include the proofs of the two combinatorial arguments
used in the main arguments and some comments on how they can be
extended.

\begin{lemm} {\rm [Proof of Lemma
\ref{justificationboxb}]}
Suppose  $\mu^{\kappa}= \mu$ and
$2^{\mu}=\mu^+ = \lambda$.
If $S^* \subseteq C^{\kappa}(\lambda)$ is stationary
then there is an $S \subseteq \lambda$ such that
$\text{\sqdi}_{\lambda,\kappa,\omega}$
holds and $S \inter C^{\kappa}(\lambda) \subseteq S^*$.
\end{lemm}

\Proof.
Let $S^* \subseteq C^{\kappa}(\lambda)$ be stationary.
Let $\langle A_{\alpha}: \alpha < \lambda\rangle$ be a list of
all bounded subsets (indeed $A_{\alpha} \subseteq \alpha$) of $\lambda$
each appearing $\lambda$ times.  (Such a list exists since
$2^{\mu} = \mu^+ = \lambda$.)
Let $\{\langle A^{\alpha,\epsilon},
C^{\alpha,\epsilon}\rangle:\epsilon < \mu\}$
list
all pairs $\langle A, C\rangle$ where $A$ has
the form $\cup_{i\in X} A_i$ with $X$ a subset of $\alpha$ of
cardinality at most $\kappa$ and $C$ is a closed subset of $\alpha$ with
$\otp(C) \leq \kappa$.   Moreover, we require that
if $\otp(C) = \kappa$ then $\alpha \in S^*$.
There are only $\mu$ such pairs since
$\mu^{\kappa} = \mu$.

By Engelking and Karlowicz \cite{EC}, there exists a sequence of
functions $F_{\xi}:\lambda \mapsto \mu$
for $\xi< \mu$ such that any partial function $h$
taking $\lambda$ to $\mu$ with $|\dom h| \leq
 \kappa$ extends to some
$F_{\xi}$.

Let $\langle\ ,\ \rangle$, $\proj_0$, $\proj_1$,  denote pairing and
projection functions on $\mu$ and on $\lambda$ such
that $\alpha$ and $\beta$ are always less than or equal to
$\langle \alpha, \beta \rangle$.

For each $\zeta < \mu$, we will define two sequences of
sets
$\Ascr^{\zeta} = \langle
A^\zeta_\alpha:\alpha < \lambda\rangle$
and $\Cscr^{\zeta}
= \langle C^\zeta_\alpha:\alpha < \lambda\rangle$ and a set $S^{\zeta}$.
We will show that
for some $\zeta$,
$\Ascr_{\zeta}$
and $\Cscr^{\zeta}$
satisfy the definition of
\boxdi$_{\lambda,\kappa,\omega}(S^{\zeta})$.

For any function $F$,
let $F^0$ and $F^1$ denote the result of applying
the first and second projection functions respectively after $F$.

If $C^{\alpha,F^0_{\zeta}(\alpha)}$ is closed and for each $\beta \in
C^{\alpha,F^0_{\zeta}(\alpha)}$, $C^{\zeta}_{\beta} = \beta \inter
C^{\alpha,F^0_{\zeta}(\alpha)}$
and, if $\alpha$ is a limit ordinal, $\alpha = \sup
C^{\alpha,F^0_{\zeta}(\alpha)}$
while each nonaccumulation point of $C^{\alpha,F^0_{\zeta}(\alpha)}$
is an even successor ordinal,
let $C^{\zeta}_{\alpha} =
C^{\alpha,F^0_{\zeta}(\alpha)}$.  Otherwise, let $C^{\zeta}_{\alpha} =
\emptyset$.

As a first approximation to $A^{\zeta}_{\alpha}$,
let $B^{\zeta}_{\alpha} = \{\gamma <
\alpha: \langle \zeta, \gamma \rangle \in
A^{\alpha,F^1_{\zeta}(\alpha)}\}$.
Now define $A^{\zeta}_{\alpha}$ by induction on $\alpha$.
At stage $\alpha$, if for each $\beta \in \acc[C^{\zeta}_{\alpha}]$,
$A^{\zeta}_{\beta} =
B^{\zeta}_{\beta}  \inter \beta$, then
$A^{\zeta}_{\alpha} =
B^{\zeta}_{\alpha}$; otherwise
$A^{\zeta}_{\alpha} =
\cup_{\beta \in \acc[C^{\zeta}_{\alpha}]} A^{\zeta}_{\beta}$.

Let $S^{\zeta}_1 = \{\delta \in S^*: \delta = \sup C^{\zeta}_{\delta}\}$.
Finally, let $S^{\zeta} = S^{\zeta}_1 \cup
\cup_{\delta \in S^{\zeta}_1}C^{\zeta}_\delta$.

Now we claim that
for some $\zeta$,
$\Ascr^{\zeta}$
and $\Cscr^{\zeta}$
satisfy the definition of
\boxdi$_{\lambda,\kappa,\omega}(S^{\zeta})$.  It is easy to check that
for each $\zeta$,
$\Ascr^{\zeta}$, $S^{\zeta}$
and $\Cscr^{\zeta}$
satisfy all the conditions except possibly ii) and vii).
(This would be possible even if all the $C^{\zeta}_{\alpha}$ were empty.)

Moreover,
$A^{\zeta}_{\alpha}
\subseteq \alpha$ and  if $\beta \in \acc[C^{\zeta}_{\alpha}]$,
$A^{\zeta}_{\beta} =
A^{\zeta}_{\alpha}  \inter \beta$.  So the nontrivial point is to
check vii) b) of \ref{boxdi}.
If for some $\zeta$
$\Ascr^{\zeta}$, $S^{\zeta}$,
and $\Cscr^{\zeta}$
work, we are finished; if not for each $\zeta$ there is
a counterexample $\langle A^*_{\zeta},E_{\zeta}\rangle$.
That is, there is some cub $D_{\zeta}$ with $D_{\zeta} \inter X_{E_{\zeta}}
= \emptyset$.  Interpreting the definition of $X_E$ (from vii) of
Definition~\ref{boxdi}),
$X_{E_{\zeta}} = \{\delta \in S^{\zeta}_1: \delta \in
E_{\zeta}\  \&\  \acc[C_{\delta}]
\subseteq E_{\zeta}\ \&
\ A^{\zeta}_{\delta} = A^*_{\zeta} \inter \delta\}$.
Since $S^{\zeta}_1 \subseteq X_{E_{\zeta}}$ showing $X_{E_{\zeta}}$
stationary also satisfies condition ii).
Let $A^* =
\{\langle \zeta,\alpha\rangle:\zeta < \mu\  \&\  \alpha \in
A^*_{\zeta}\}$. (Note that by the definition of $\langle\ ,\ \rangle$,
this is a set of ordinals).  Let $D = \inter_{\zeta < \mu}D_{\zeta}$.

Let $E$ be the set of $\delta \in \inter_{\zeta<\mu}E_{\zeta}$
 such that $\delta$ is closed
under the pairing and projection functions and for each $\alpha <
\delta$, $A^* \inter \alpha \in \{A_i:i< \delta\}$.

Define a function $g:\lambda \mapsto \lambda$ by $g(i)$ is the least $j$
such that $A^* \inter i = A_j$.  Since the $A_i$ enumerate all bounded
subsets of $\lambda$ this function is well defined.
Note that $i < \alpha$ is in $E$ implies $g(i) < \alpha$.

Then $E$ is a cub
in $\lambda$.
To demonstrate a contradiction, fix
$\delta \in S^* \inter C^{\kappa}(\lambda) \inter \acc[E] \inter D$
($\delta$ exists
as $S^* \subseteq C^{\kappa}(\lambda)$ is stationary)
and then choose a subset $C_{\delta}$ of $\delta$ such that
$\acc [C_{\delta}] \subseteq E$,
each nonaccumulation point of
$C_{\delta}$ is an even successor ordinal, and $\delta = \sup C_{\delta}$.

Let $Y_{\delta}$ denote
$\acc[C_{\delta}] \union \{\delta\}$.

For each $\alpha \in Y_{\delta}$, $X_{\alpha} = \{g(i): i \in C_{\delta}
\inter \alpha \}$ is a subset of $\alpha$ with cardinality less than
$\kappa$ and so for some $\nu$, $A^{\alpha,\nu} = \union_{j
\in X_{\alpha}}
 A_j$.  But,
$$\union_{j\in X_{\alpha}}A_j
= \union_
{i \in C_{\delta}\inter
\alpha}C_{g(i)}
= \union_
{i \in C_{\delta}\inter
\alpha}
(A^* \inter i)=
(A^* \inter \alpha)$$
(as $\alpha \in Y_{\delta}$); so
$A^* \inter \alpha =
A^{\alpha,\nu}$.
By the choice of $\langle C^{\alpha,\epsilon}:\epsilon < \mu\rangle$,
for some $\nu'$, $C_{\delta} \inter \alpha =
C^{\alpha,\nu'}$.

Based on these two observations we can now define a function $h$
with domain $\{\delta\} \union C_{\delta}$ as follows.
For each $\gamma \in Y_{\delta}$,
let $h(\gamma)$ be the minimal
$\epsilon < \mu$ with
$A^* \inter \gamma =
A^{\gamma,\proj_1(\epsilon)}$ and $C_{\delta} \inter \gamma =
C^{\gamma,\proj_0{(\epsilon)}}$.
For $\gamma \in \nacc[C_{\delta}]$, we only care about $h_0(\gamma)$;
let $h(\gamma)$ be the minimal
$\epsilon < \mu$ with
$C_{\delta} \inter \gamma =
C^{\gamma,\proj_0(\epsilon)}$.

Now for some $\zeta$, the choice of the $F_{\zeta}$ guarantees that $h
\subseteq F_{\zeta}$.  For this $\zeta$,
we show by induction on $\alpha \in C_{\delta}\union\{\delta\}$ that
$C^{\zeta}_{\alpha} = C^{\alpha,h^0(\alpha)}$
(and thus $\delta \in S^{\zeta}_1$)
and $A^{\zeta}_{\delta} = A^*_{\zeta}
\inter \delta$ (so $\delta \in X_{E_{\zeta}}$).
This contradicts the choice of $D_{\zeta}$ and completes the proof.

We first show $C^{\zeta}_{\alpha} = C^{\alpha,h^0(\alpha)}$ by induction on
$\alpha \in C_{\delta}\union \{\delta\}$.
For each such $\alpha$, $C_{\delta} \inter \alpha = C^{\alpha,h^0(\alpha)}$.
If $\alpha$
is least in $C_{\delta}$, this implies $C^{\zeta}_{\alpha} =
C^{\alpha,h^0(\alpha)}$.
Suppose the result holds for $\beta \in C_{\delta}$ that are less
than $\alpha$.  Then by induction for $\beta < \alpha$,
$$C^{\zeta}_{\beta}
=C^{\beta,h^0(\beta)}
= C_{\delta} \inter \beta
= (C_{\delta} \inter \alpha) \inter \beta
=C^{\alpha,h^0(\alpha)} \inter \beta.$$
The requirements of the definition of $C^{\zeta}_{\alpha}$ are met
so $C^{\zeta}_{\alpha} = C^{\alpha,h^0(\alpha)}$ as required.

A similar induction shows
$A^{\zeta}_{\alpha} = A^*_{\zeta} \inter \alpha$.  Suppose first that
$\alpha$ is minimal in $C_{\delta}$.   Then by definition,
$A^{\zeta}_{\alpha} = \{\gamma:\langle \zeta, \gamma\rangle \in
A^{\alpha, F^1_{\zeta}(\alpha)}\}$.  But $h \subseteq F_{\zeta}$
so $A^* \inter \alpha = A^{\alpha, F^1_{\zeta}(\alpha)}$.  Therefore,
$$A^{\zeta}_{\alpha} =
B^{\zeta}_{\alpha} = \{\gamma: \langle \zeta,\gamma\rangle \in A^*
\inter \alpha \} = A^*_{\zeta} \inter \alpha.$$  Now suppose that for
$\beta \in C_{\delta}$ with $\beta < \alpha$,
$$A^{\zeta}_{\beta} = A^*_{\zeta} \inter \beta = B^{\zeta}_{\beta}.$$
To see this equality holds for $\alpha$ note that by definition of
$A^{\zeta}_{\alpha}$, this yields $A^{\zeta}_{\alpha} =
B^{\zeta}_{\alpha}$.  But $B^{\zeta}_{\alpha} = A^*_{\zeta} \inter
\alpha$ since $h \subseteq F_{\zeta}$.

We have established Lemma~\ref{justificationboxb} in ZFC.  The
specific requirements on $\kappa$ and $\lambda$ can be weakened
if $V = L$.
We would like to remove the restriction that $R = \omega$ and that
$\mu^{\kappa} = \mu$.  The principal {\bf SD}$_{\lambda}$ is
established for all uncountable successor cardinals $\lambda
= \mu^+$ in
\cite{Shelah221} assuming $V = L$.
This principal differs from
\boxdi$_{\lambda,\kappa,R}(S)$ in three ways:
the cubs $C_{\alpha}$ are defined only for limit ordinals $\alpha$,
the cardinal $\kappa$ is identified with $\mu$, there is no explicit
treatment of $R$.  In fact, {\bf SD}$_{\lambda}$ describes the default case
where $R = \lambda$.  Since we assume $\mu^{\kappa} = \mu$, our proof of
Lemma~\ref{justificationboxb} definitely misses the case $\kappa = \mu$.

We outline
technical modifications of a system $\Cscr$ satisfying {\bf SD}$_\lambda$
to show

\begin{lemm}[V=L]
{\rm [Proof of Lemma
\ref{justificationboxc}]}
If $\lambda = \mu^+$ and $\kappa \leq \mu$, then
for some stationary $S$,
\boxdi$_{\lambda,\kappa,\lambda}(S)$.
\end{lemm}

\Proof.  We show how to use to $C_{\alpha}$ from {\bf SD}$_{\lambda}$
to construct first a sequence $\Cscr^1$ such that
\boxdi$_{\lambda,\mu,\lambda}(S)$ holds with $S = \lambda$.
Let $E$ be the set of ordinals $> 0$ in $\lambda$ divisible by $\mu$
and let $E^*$ denote the accumulations points of $E$.

We define by induction on $\alpha \in E$
a sequence, increasing with $\alpha$, $\Cun^1|(\alpha +1) =
\langle C^1_{\beta}:\beta \leq \alpha\rangle$.

\begin{enumerate}
\item $\alpha = \min E$:  Thus, $\alpha = \mu$.  If $\beta$ is even
 $C^1_{\beta}$ is the
 set of even ordinals less than $\beta$; if $\beta$ is
odd,  $C^1_{\beta}$ is empty.
\item $\alpha$ is a successor in $E$:
Thus $\alpha = \alpha^* + \chi$
where $\alpha^* \in E$ and $\Cun^1|(\alpha^* +1)$ has been defined.
Fix a map $h_{\alpha}$ from $\mu$ onto $\alpha^*$.
\begin{enumerate}
\item $\alpha^* < \beta \leq \alpha$ and $4|\beta$:
$C^1_{\beta}$ consists of those ordinals in $(\alpha^*,\beta)$ that
are divisible by $4$.
\item $\alpha^* < \beta \leq \alpha$ and $\beta$ is odd:
$C^1_{\beta}$ is empty.
\item $\alpha^* < \beta \leq \alpha$ and $4\not|\beta$ but $2|\beta$:
Thus $\beta$ has the form $\alpha^* + 4i + 2$ for some $i$.
\begin{itemize}
\item $\otp C^1_{h(i)} < \mu$ and $h(i)$ is even:
$C^1_{\beta} = C^1_{h_1(i)} \union \{h(i)\}$
\item otherwise: $C^1_{\beta}$ is empty.
\end{itemize}
\end{enumerate}
\item $\alpha = \delta\in E^*$, i.e a limit in $E$: we only
have to define $C^1_{\delta}$.
\begin{enumerate}
\item $\acc[C_{\delta}] \inter E^* = \emptyset$:  Necessarily $\cf(\delta) =
\aleph_0$.  Choose an increasing
sequence
$\{\alpha_n:n <\omega\}$ from $E$ with limit $\alpha$.
Then choose by induction on $n$ successor even ordinals $\beta_n$ with
$\alpha_n < \beta_n < \alpha_n + \mu$ such that $C^1_{\beta_m} = \{\beta_n:
n < m\}$.  (The third clause of the previous case is the key to this
induction.)
Finally, let $C^1_{\delta}$ be the set of $\beta_n$.
\item $\acc[C_{\delta}] \inter E^* \neq \emptyset$ but $\acc[C_{\delta}]
\inter E^*$ is bounded in $\delta$ by some $\delta'$.  Choose $\alpha_n$
with limit $\delta'$ and $\beta_n$ as before but so that
$C^1_{\beta_n} = C_{\delta'} \union \{\delta'\}\union \{\beta_m:m <n\}$.
Then let
$C^1_{\delta} = C_{\delta'} \union \{\delta'\}\union \{\beta_m:m <
\omega\}$.
\item $\delta = \sup( \acc[C_\delta] \inter E^*)$: Let $C^1_{\delta} =
\union_{\alpha \in C_{\delta} \inter E} C^1_{\alpha}$.

\end{enumerate}
\end{enumerate}

Now one shows by induction on $\alpha$ that if $\beta \in
 C^1_{\alpha}$, then $C^1_{\alpha} \inter \beta = C^1_{\beta}$.
(The other requirements on the $C^1_{\alpha}$ are easily verified.)

If $\kappa < \mu$ we modify the $C^1_{\alpha}$ defined in the first stage
of the proof as follows.  If $\otp C^1_{\alpha} \leq \kappa$ then
$C^2_{\alpha} = C^1_{\alpha}$.
If $\otp C^1_{\alpha} > \kappa$ then
$C^2_{\alpha} = \{\beta \in C^1_{\alpha}:\otp C^1_{\beta} > \kappa\}$.
Let $S_1 = \{\alpha < \lambda: \cf(\alpha) = \kappa \ \& \ \otp C^1_{\alpha}
= \kappa\}$.
And $\tilde S_1 = \{\alpha < \lambda: \otp C^1_{\alpha} \leq \kappa\}$.

\jtbnumpar{Concluding Remarks}
\mbox{}
\begin{enumerate}
\item
Examination of our proof
of Lemma~\ref{justificationboxb}
shows that in fact given a square sequence, it is
possible to add on a diamond sequence
to satisfy \boxdi$_{\lambda,\kappa,R}(S)$.
We state this explicitly in Lemma~\ref{boxextend} below.
\item The proof of Lemma~\ref{justificationboxc} can be extended
to allow for $R < \lambda$.
\item In fact, the derivation of
$\text{\sqdi}_{\lambda,\kappa,R}$
for $R \leq \lambda$ from the assumption that $\mu = \mu^{\kappa}$,
$\lambda = 2^{\mu}$ and there is a square on
$\{\delta<\lambda:\cf(\delta) < R\}$ but
without assuming V=L will be published
elsewhere by the second author.
\item Similar results hold assuming V=L for $\lambda$ inaccessible and
$\lambda > \kappa$ with $R \leq \lambda$ should follow by the methods of
Beller and Litman \cite{BellerLitman} but we have not checked this in detail.
\end{enumerate}

\begin{lemm}
\label{boxextend}
Suppose
$2^{\mu} = \mu^+ = \lambda$, $\mu^{\kappa} = \mu$,
and $S$, $R$, $\langle C_{\delta}:\delta \in S\rangle$ have been chosen
to satisfy all conditions of
\boxdi$_{\lambda,\kappa,R}(S)$ except ii)
and vii).
Suppose further that for each cub $E$ of $\lambda$,
$\{\delta: \acc[C_{\delta}] \subseteq E\}$is stationary in $\lambda$.
Then \boxdi$_{\lambda,\kappa,R}(S)$ holds.
\end{lemm}
\bibliography{ssgroups}
\bibliographystyle{plain}
\end{document}